\input amstex
\documentstyle{amsppt}
\magnification 1200

\def\diam{\operatorname{diam}}

\def\N{\Bbb N}
\def\Z{\Bbb Z}

\def\S{\Cal S}
\def\ex{\operatorname{ex}}
\def\cl{\operatorname{cl}}
\def\dim{\operatorname{dim}}

\def\Ind{\operatorname{Ind}}
\def\asdim{\operatorname{asdim}}
\def\asInd{\operatorname{asInd}}
\def\trInd{\operatorname{trInd}}
\def\trasInd{\operatorname{trasInd}}
\topmatter
\title
Addition and subspace theorems for asymptotic large inductive
dimension
\endtitle
\rightheadtext {Addition and subspace theorems}
\author
T.Radul
\endauthor
\address
Dept. de Matematicas, Facultad de Cs. Fisicas y Mat., Universidad
de Concepcion, CASILLA 160-C, Concepcion, Chile
 \newline
e-mail: tarasradul\@yahoo.co.uk
\endaddress

\keywords Asymptotic dimension,transfinite extension
\endkeywords
\subjclass 54F45, 54D35
\endsubjclass
\abstract We prove addition and subspace theorems for asymptotic
large inductive dimension. We investigate a transfinite extension
of this dimension and show that it is trivial.
\endabstract
\endtopmatter
\document
\baselineskip18pt

{\bf 0.} Asymptotic dimension $\asdim$ of a metric space was
defined by Gromov for studying asymptotic invariants of discrete
groups [1]. This dimension can be considered as asymptotic
analogue of the Lebesgue covering dimension $\dim$. Dranishnikov
has introduced dimension $\asInd$ which is analogous to large
inductive dimension $\Ind$ [2].

Some of  basic theorems of classical dimension theory are sum,
addition and subspace theorems for different dimensions and
classes of topological spaces. Here we mention some of them
related to dimension $\Ind$. (All mentioned facts from classic
dimension theory could be found in [3]).

1. Countable sum theorem: If a strongly hereditarily normal space
$X$ can be represented as the union of a sequence $F_1, F_2,\dots$
of closed subsets such that $\Ind F_i\le n$ for each $i\in \N$,
then $\Ind X\le n$.

2. Addition theorem: If a hereditarily normal space $X$ is
represented as the union of two subspaces $X_1$ and $X_2$, then
$\Ind X\le \Ind X_1+\Ind X_2+1$.

Let us remark that we cannot generalize the above theorems for the
class of all normal spaces. We have only a weaker result:

3. If a normal space $X$ is represented as the union of two closed
subsets $X_1$ and $X_2$, then $\Ind X\le \Ind X_1+\Ind X_2$.

And finally:

4. Subspace theorem: For each closed subset $M$ of a normal space
$X$ we have $\Ind M\le \Ind X$. Let us remark that if $X$ is a
strongly hereditarily normal space, then the condition that $M$ is
closed could be dropped.

There is no countable sum theorem for $\asInd$. Indeed, the space
of enters $\Z$ has asymptotic dimension $1$ but it is  countable
union of its points which have asymptotic dimension $-1$. However
we prove subspace and addition theorem for $\asInd$ in this paper.

Extending codomain of $\Ind$ to ordinal numbers we obtain the
transfinite extension $\trInd$ of the dimension $\Ind$. It is
known that there exists a space $S_\alpha$ such that $\trInd
S_\alpha=\alpha$ for each countable ordinal number $\alpha$.
Zarichnyi has proposed to consider transfinite extension of
$\asInd$ and conjectured that this extension is trivial. We prove
this conjecture: if a space has a transfinite asymptotic
dimension, then its dimension is finite.

The paper is organized as follows: in Section 1 we give some
necessary definitions and introduce some denotations, in Section 2
we prove some theorem which can be consider as a weak version of
countable sum theorem, in Section 3 we obtain the main results and
in Section 4 we show that the transfinite extension of $\asInd$ is
trivial.

{\bf 1.} Let $A_1,A_2\subset X$ be two disjoint closed subsets in
a topological space. We recall that a {\it partition} between
$A_1$ and $A_2$ is a subset $C\subset X$ such that there are open
disjoint sets $U_1,U_2$ satisfying the conditions: $X\setminus
C=U_1\cup U_2$, $A_1\subset U_1$ and $A_2\subset U_2$. Clearly a
partition $C$ is a closed subset of $X$.

We recall the definition of the large inductive dimension $\Ind$
[3]: $\Ind X=-1$ iff $X=\emptyset$; $\Ind X\le n$ if for every two
disjoint closed subsets $A_1, A_2\subset X$ there is a partition
$C$ with $\Ind C\le n-1$.

We will define the dimension $\asInd$ for the class of proper
metric space. We recall that a metric space is {\it proper}  if
every closed ball is compact. We assume that there is chosen some
base point $x_0\in X$ for each proper metric space $X$. The
generic metric we denote by $d$. If $X$ is a metric space and
$A\subset X$ we denote by $N_r(A)$ the closed $r$-neighborhood of
$A$ in $X$: $N_r(A)=\{x\in X\mid d(x,A)\le r\}$ and by $B_r(A)$
the open $r$-neighborhood: $B_r(A)=\{x\in X\mid d(x,A)< r\}$.

A subset $W$ of a metric space is called an {\it asymptotic
neighborhood} of a set $A\subset X$ if
$\lim_{r\to\infty}d(A\setminus B_r(x_0),X\setminus W)=\infty$. We
call two subsets $A_1, A_2\subset X$ in a metric space $X$ {\it
asymptotically disjoint} if $\lim_{r\to\infty}d(A_1\setminus
B_r(x_0),A_2\setminus B_r(x_0))=\infty$. It is easy to see that
sets $A_1$ and $A_2$ are asymptotically disjoint iff $X\setminus
A_2$ is an asymptotic neighborhood of $A_1$ and $X\setminus A_1$
is an asymptotic neighborhood of $A_2$.

A map $\phi:X\to I=[0,1]$ is called {\it slowly oscillating} if
for any $r>0$, for given $\varepsilon>0$ there exists $D>0$ such
that $\diam \phi(B_r(x))<\varepsilon$ for any $x$ with
$d(x,x_0)\ge D$. If $C_h(X)$ is the set of all continuous slow
oscillating functions $\phi:X\to I$, then the {\it Higson
compactification} is the closure of the image of $X$ under the
embedding $\Phi:X\to I^{C_h(X)}$ defined as
$\Phi(x)=(\phi(x)\mid\phi\in C_h(X))\in I^{C_h(X)}$. We denote the
Higson compactification of a proper metric space $X$ by $cX$ and
the remainder $cX\setminus X$ by $\nu X$. The compactum $\nu X$ is
called {\it Higson corona}. Let us remark that $\nu X$ does not
need to be metrizable.

Let $C$ be a subset of a proper metric space $X$. By $C'$ we
denote the intersection $\cl C\cap \nu X$ of the closure $\cl C$
in the Higson compactification $cX$. Clearly, two sets $A_1$ and
$A_2$ are asymptotically disjoint iff their traces $A_1'$ and
$A_2'$ in the Higson corona are disjoint. Note that for each $r>0$
we have $N_r(C)'=C'$.

Let $A_1, A_2\subset X$ be two asymptotically disjoint subsets of
a proper metric space $X$. A closed subset $C\subset X$ is called
an {\it asymptotic separator} for $A_1$ and $A_2$ if its trace
$C'$ is a partition for $A_1'$ and $A_2'$ in $\nu X$.

We define $\asInd X=-1$ if and only if $X$ is bounded; $\asInd
X\le n$ if for every two asymptotically disjoint sets $A$,
$B\subset X$ there is an asymptotic separator $C$ with $\asInd
C\le n-1$. Naturally we say $\asInd X=n$ if $\asInd X\le n$ and it
is not true that $\asInd X\le n-1$. We set $\asInd X=\infty$ if
$\asInd X>n$ for each $n\in \N$ [2].

For each $A\subset Y\subset X$ we denote by $\ex_Y
A=Y'\setminus(Y\setminus A)'$. Clearly, $\ex_Y A$ is an open set
in $Y'$. By $\omega^\omega$ we denote the set of all functions
$\tau:\N\cup\{0\}\to\N\cup\{0\}$ such that $\tau(0)=0$. For each
$\tau\in\omega^\omega$ we define the set $V_\tau^Y(A)$ as follows:
$V_\tau^Y(A)=\{y\in Y\mid d(y,x_0)\ge\tau([d(y,A)])\}$ where
$[d(y,A)]$ is the enter part of real number $d(y,A)$. Clearly,
$A\subset V_\tau^Y(A)$ for each $\tau\in\omega^\omega$. If $X=Y$
we use simpler denotations $\ex A$ and $V_\tau(A)$.

\proclaim {Lemma 1} The family
$\{\ex(V_\tau(A))\mid\tau\in\omega^\omega\}$ forms the base of
neighborhoods of the set $A'$ in the space $\nu X$.
\endproclaim

\demo{Proof} Let us show firstly that $A'\subset \ex(V_\tau(A))$
for each $\tau\in\omega^\omega$. It is enough to show that $A$ and
$X\setminus V_\tau(A)$ are asymptotically disjoint. Fix any $D>0$.
Put $R=\max\{\tau(l)\mid l\in\{0,\dots,[D]\}\}$. Then for each
$x\in (X\setminus V_\tau(A))\setminus B_R(x_0)$ we have $R\le
d(x,x_0)<\tau([d(x,A)])$. Then $d(x,A)\ge [D]+1>D$.

Consider now any closed subset $B$ of $\nu X$ such that $B\cap
A'=\emptyset$. Then there exists a continuous function $f:\nu X\to
[0,1]$ such that $f(B)\subset \{0\}$ and $f(A')\subset \{1\}$. We
can extend $f$ to a continuous function $g:cX\to [0,1]$ such that
$A\subset g^{-1}(1)$. Put $C=g^{-1}[0,\frac{1}{2}]\cap X$. The
sets $C$ and $A$ are asymptotically disjoint. For each $n\in\N$
there exists $R(n)>0$ such that $d(C\setminus B_{R(n)}(x_0),A)\ge
n$. Put $\tau(0)=0$ and $\tau(n)=[R(n+1)]$. Choose any $c\in C$.
Then $d(c,x_0)<R([d(c,A)]+1)=\tau([d(c,A)])$. So, $C\subset
X\setminus V_\tau(A)$ and $B\subset C'\subset (X\setminus
V_\tau(A))'\subset \nu X\setminus\ex V_\tau(A)$. Hence
$\{\ex(V_\tau(A))\mid\tau\in\omega^\omega\}$ forms the base of
neighborhoods of the set $A'$ in the space $\nu X$ and the lemma
is proved.
\enddemo

{\bf 2.} Let $X$ be a proper metric space and $X_0$ is an
unbounded subset of  $X$. We say that $X_0$ is a {\it kernel} of
$X$ if there exists a sequence $(k_i)_{i=0}^\infty$ of natural
numbers such that $k_i\to\infty$ and for each $i\in\N\cup\{0\}$,
$x\in X\setminus N_i(X_0)$ we have $B_{k_i}(x)=\{x\}$. We suppose
that $x_0\in X_0$ where $x_0$ is the base point in $X$.

\proclaim {Lemma 2} If $X_0$ is a kernel of a metric proper space
$X$ then the family $\{V_\tau(X_0))\mid\tau\in\omega^\omega\}$
forms the base of clopen neighborhoods of the set $X_0'$ in the
space $\nu X$.
\endproclaim

\demo {Proof} It follows from Lemma 1 that it is enough to prove
that $V_\tau(X_0)'\cap (X\setminus V_\tau(X_0))'=\emptyset$.
Suppose the contrary: there exist $x\in V_\tau(X_0)'\cap
(X\setminus V_\tau(X_0))'$. Let $U$ be a neighborhood of $x$ in
$cX$. Then there exist two sequences $(a_i)$ and $(b_i)$ in
$V_\tau(X_0)\cap U$ and $(X\setminus V_\tau(X_0))\cap U$
respectively such that $0<d(a_i,b_i)\le r$ for some $r>0$ and
$a_i$, $b_i\in X\setminus B_i(x_0)$. Choose any $n_0\in \N$ such
that $k_n>r$ for each $n\ge n_0$ where $(k_n)$ is a sequence from
the definition of kernel. Then $a_i$, $b_i\in N_{n_0}(X_0)$ for
each $i\in \N$. Hence $\emptyset\neq \cl U\cap N_{n_0}(X_0)'=\cl
U\cap X_0'$ and $x\in X_0$. We obtain the contradiction.
\enddemo

\proclaim {Lemma 3} If $X_0$ is a kernel of $X$ then $\Ind (\nu
X\setminus V_\tau(X_0)')\le 0$ for each $\tau\in\omega^\omega$.
\endproclaim

\demo {Proof} Since $\nu X\setminus V_\tau(X_0)'$ is compact, it
is enough to prove that the space $\nu X\setminus V_\tau(X_0)'$
has a base of clopen sets [3, Theorem 1.6.5].

Choose any $x\in\nu X\setminus V_\tau(X_0)'$ and its open
neighborhood $U\subset \nu X\setminus V_\tau(X_0)'$. Take a
continuous function $f:\nu X\to[0,1]$ such that $f(x)=0$ and
$f(\nu X\setminus U)\subset\{1\}$. Extend $f$ to the continuous
function $g:cX\to [0,1]$. Put $A=g^{-1}[0,\frac{1}{3}]\cap X$ and
$C=g^{-1}[\frac{2}{3},1]\cap X$. The sets $C$ and $A$ are
asymptotically disjoint.  Moreover, $x\in A'$ and $\nu X\setminus
U\subset C'$. For each $n\in\N$ there exists $R(n)>0$ such that
$d(C\setminus B_R(x_0),A\setminus B_R(x_0))\ge n$. Put $\tau(0)=0$
and $\tau(n)=[R(n+1)]$. We can show that $V_\tau(A)'$ is a clopen
neighborhood of $x$ such that $V_\tau(A)'\subset U$ using the same
reasoning as in Lemmas 1 and 2. The lemma is proved.
\enddemo

Let us define a preorder $\le^*$ in $\omega^\omega$ as follows
$\tau\le^*\sigma$ iff there exists $n\in\N$ such that
$\tau(i)\le\sigma(i)$ for each $i\ge n$. It is easy to check that
$V_\tau(A)'\subset V_\sigma(A)'$ if $\sigma\le^*\tau$.

\proclaim {Theorem 1} Let  $X_0$ be a kernel of a metric proper
space $X$ such that $\asInd X_0\le k\ge 0$. Then $\asInd X\le k$.
\endproclaim

\demo {Proof} We use  induction respect to $k$. Let $\asInd X_0\le
0$. Then $\Ind X_0'\le 0$ [2]. Consider any two asymptotically
disjoint sets $A$ and $B$ in $X$. We can represent $X_0'=K\cup L$
where $K$ and $L$ are two disjoint closed subsets of $X_0'$ such
that $A'\cap X_0'\subset K$ and $B'\cap X_0'\subset L$. Then
$A'\cup K$ and $B'\cup L$ are two disjoint closed subsets of $\nu
X$. Choose two disjoint open subsets $U_1$ and $U_2$ of $\nu X$
such that $A'\cup K\subset U_1$ and $B'\cup L\subset U_2$. Then
$U_1\cup U_2$ is a neighborhood of $X_0$ in $\nu X$ and there
exists $\tau\in\omega^\omega$ such that $V_\tau(X_0)' \subset
U_1\cup U_2$. The sets $A'\cap (\nu X\setminus V_\tau(X_0)')$ and
$B'\cap (\nu X\setminus V_\tau(X_0)')$ are two disjoint closed
subset of the 0-dimensional  space $\nu X\setminus V_\tau(X_0)'$,
so, there exist two open disjoint subsets $O_1$, $O_2$ of $\nu
X\setminus V_\tau(X_0)'$ such that $A'\cap (\nu X\setminus
V_\tau(X_0)')\subset O_1$, $B'\cap (\nu X\setminus
V_\tau(X_0)')\subset O_2$ and $O_1\cup O_2=\nu X\setminus
V_\tau(X_0)'$. Put $V_1=(U_1\cap V_\tau(X_0)')\cup O_1$ and
$V_2=(U_2\cap V_\tau(X_0)')\cup O_2$. Since $V_\tau(X_0)'$ is a
clopen subset of $\nu X$, the sets $V_1$ and $V_2$ are open.
Moreover, the sets $V_1$ and $V_2$ are disjoint, $V_1\cup V_2=\nu
X$ and $A'\subset V_1$ and $B'\subset V_2$. Thus, empty space is a
partition between $A'$ and $B'$ in $\nu X$. Hence, empty space is
an asymptotic separator between $A$ and $B$ in $X$ and $\asInd
X\le 0$.

Suppose that the theorem is proved for each $i<n\ge 1$. Consider
the case when $\asInd X_0\le n$. Let $A$ and $B$ be any
asymptotically disjoint subset of $X$. Then $A'\cap X_0'$ and
$B'\cap X_0'$ are disjoint closed subsets of $X_0'$ and we can
choose a continuous function $f': X_0'\to [0,1]$ such that
$f'(A'\cap X_0')\subset \{0\}$ and $f'(B'\cap X_0')\subset \{1\}$.
We can extend this function to a continuous function $f:\cl X_0\to
[0,1]$. The sets $A_1=(f^{-1}[0,\frac{1}{3}])\cap X_0$ and
$B_1=(f^{-1}[\frac{2}{3},1])\cap X_0$ are asymptotically disjoint
and we can choose an asymptotic separator $L\subset X_0$ between
them such that $\asInd L<n$.

Represent $X=\cup_{i=0}^\infty X_i$ where $X_i=N_i(X_0)\setminus
N_{i-1}(X_0)$ for $i\in\N$. It follows from the definition of
kernel that for each $R>0$ there exists $i(R)\in \N$ such that the
set $\cup_{k=i(R)}^\infty X_k$ is $R$-discrete and
$d(\cup_{k=i(R)}^\infty X_k,\cup_{k=0}^{i(R)-1} X_k)\ge R$.

Since $L$ is an asymptotic separator in $X_0$ between $A_1$ and
$B_1$, then $L'$  is a partition in $X_0'$ between $A_1'$ and
$B_1'$. Then we can choose two open disjoint sets $O_A$ and $O_B$
in $X_0'$ such that $A_1'\subset O_A$, $B_1'\subset O_B$ and
$X_0'\setminus L'= O_A\cup O_B$.

For each $\tau\in\omega^\omega$ we can represent $X_0\setminus
V_\tau(L)$ as a union of two disjoint sets $A_\tau$ and $B_\tau$
such that $A_\tau'\subset O_A$ and $B_\tau'\subset O_B$. Moreover,
we can suppose that for each $\tau\le^*\sigma$ there exists $R>0$
such that $A_\sigma\supset A_\tau\setminus B_R(x_0)$ and
$B_\sigma\supset B_\tau\setminus B_R(x_0)$. Define for each
$\tau\in\omega^\omega$ two subsets $C_\tau$, $D_\tau$ in $X$ as
follows: $C_\tau=\cup_{i=0}^\infty\{x\in X_i\mid N_i(x)\cap
X_0\subset A_\tau\}$ and $D_\tau=\cup_{i=0}^\infty\{x\in X_i\mid
N_i(x)\cap B_\tau\neq\emptyset\}$. We have $C_\tau\cap
D_\tau=\emptyset$ and $X\setminus (C_\tau\cup D_\tau)\subset
\cup_{i=0}^\infty\{x\in X_i\mid N_i(x)\cap
V_\tau^{X_0}(L)\neq\emptyset\}$.

Let us show that $\ex C_\tau\supset\ex_{X_0}A_\tau$. Choose any
point $x\in\ex_{X_0}A_\tau$. Then there exists  $Z\subset X_0$
such that $Z$ and $X_0\setminus A_\tau$ are asymptotically
disjoint and $x\in Z'$. Choose any $a>0$. Since $k_i\to\infty$,
there exists $n_0\in\N$ such that $k_{n_0}\ge a$. Choose any $R>0$
such that $d(Z\setminus B_R(x_0),(X_0\setminus A_\tau)\setminus
B_R(x_0))\ge a+n_0$. Consider any points $z\in Z\setminus
B_{R+n_0}(x_0)$ and $y\in\cup_{i=0}^\infty\{x\in X_i\mid
N_i(x)\cap(X_0\setminus A_\tau) \neq\emptyset\}$. We have that
$d(z,y)\ge a$. Hence the sets $Z$ and $\cup_{i=0}^\infty\{x\in
X_i\mid N_i(x)\cap(X_0\setminus A_\tau) \neq\emptyset\}$ are
asymptotically disjoint and $x\in\ex C_\tau$. Analogously we can
show that $\ex D_\tau\supset\ex_{X_0}B_\tau$.

Now consider any $x\in C_\tau'\cap D_\tau'$. Then for every
neighborhood $V$ of $x$ in the Higson compactification $cX$ there
exists two sequences $(c_i)$ in $V\cap C_\tau$ and $(d_i)$ in
$V\cap D_\tau$ such that $0<d(c_i,d_i)\le r$ for some $r>0$ and
$c_i$, $d_i\in X\setminus B_i(x_0)$. Choose $n_0\in\N$ such that
$k_{n_0+1}>r$. So, $c_i$, $d_i\in\cup_{k=0}^{n_0} X_k$ and
$\emptyset\neq\cl V\cap(\cup_{k=0}^{n_0} X_k)'=\cl V\cap X_0'$.
Hence we have $x\in X_0'$. Moreover, $x\in X_0'\setminus
(\ex_{X_0}A_\tau\cup\ex_{X_0}B_\tau)=V_\tau^{X_0}(L)'$.

Denote $S_k(M)=\cup_{i=0}^\infty\{x\in X_i\mid N_{ki}(x)\cap M
\neq\emptyset\}$ for any $M\subset X_0$ and $k\in\N$. The set
$S_1(M)$ we denote simply by $S(M)$. We have that $\nu X\setminus
(\ex D_\tau\cup \ex C_\tau)=(X\setminus D_\tau)'\cap(X\setminus
C_\tau)'\subset(X\setminus(D_\tau\cup C_\tau))'\cup(D_\tau'\cap
C_\tau')\subset (S(V_\tau^{X_0}(L)))'$.

Put $K_j=S(B_j(x_0))$ for each $j\in\N$. We have that the sets
$K_j$ and $X_0$ are asymptotically disjoint for each $j\in \N$.
There exists $\sigma_j\in\omega^\omega$ such that
$V_{\sigma_j}(X_0)'\cap K_j'=\emptyset$. Define
$\sigma_0\in\omega^\omega$ as follows $\sigma_0(0)=0$ and
$\sigma_0(i)=\max\{\sigma_j(i)\mid j\le i\}$. We have that
$\sigma_j\le^*\sigma_0$ for each $j\in\N$. Hence
$\cl\cup_{j=1}^\infty K_j'\subset\nu X\setminus
V_{\sigma_0}(X_0)'$.

Let us show that
$\cap_{\tau\in\omega^\omega}(S(V_\tau^{X_0}(L)))'\cap
V_{\sigma_0}(X_0)'\subset (S_2(L))'$. Choose any $x\notin
(S_2(L))'\cup(\nu X\setminus V_{\sigma_0}(X_0)')$. Then there
exists a subset $Z\subset X$ such that $x\in Z'$ and the sets $Z$,
$S_2(L)\cup(X\setminus V_{\sigma_0}(X_0))$ are asymptotically
disjoint. Let us consider the set $O(Z)\subset X_0$ defined as
follows $O(Z)=\{y\in X_0\mid$ there exists $i\in\N\cup\{0\}$ and
$x\in X_i\cap Z$ such that $d(y,x)\le i\}$.

Suppose that there exist $r\in\N$ and two sequences $(y_i)$ in
$O(Z)$ and $(l_i)$ in $L$ such that $d(y_i,l_i)\le r$ and $l_i$,
$y_i\in X_0\setminus B_i(x_0)$. For each $i$ choose $x_i^{j_i}\in
X_{j_i}\cap Z$ such that $d(x_i^{j_i}, y_i)\le j_i$. Consider two
cases:

1. there exists $n_0\in \N$ such that $j_i\le n_0$ for each
$i\in\N$. Then $d(x_i^{j_i},l_i)\le n_0+r$ and we obtain the
contradiction with asymptotic disjointness of the sets $Z$ and
$S_2(L)$.

2. in the contrary case we can suppose that $j_i\ge r$ for each
$i\in \N$ and $j_i\to\infty$. Then we have $d(x_i^{j_i},l_i)\le
j_i+r\le 2j_i$. Since $x_i^{j_i}\in X_{j_i}$, we have
$x_i^{j_i}\in S_2(L)\cap Z$ and we obtain the contradiction again.
So, the sets $O(Z)$ and $L$ are asymptotically disjoint.

We can choose some $\tau\in\omega^\omega$ such that $O(Z)$ and
$V_\tau^{X_0}(L)$ are asymptotically disjoint. Let us show that
$Z$ and $S(V_\tau^{X_0}(L))$ are asymptotically disjoint. Suppose
the contrary. Then there exist $r\in\N$ and two sequences $(z_i)$
in $Z$ and $(s_i)$ in $S(V_\tau^{X_0}(L))$ such that
$d(z_i,s_i)\le r$ and $z_i$, $s_i\in X_0\setminus B_i(x_0)$.
Consider two cases:

1. there exists $n_0\in\N$ such that $z_i$,
$s_i\in\cup_{i=0}^{n_0} X_i$. We can choose $y_i\in O(Z)$, $l_i\in
V_\tau^{X_0}(L)$ for each $i\in\N$ such that $d(y_i,z_i)\le n_0$
and $d(l_i,s_i)\le n_0$. Then we have $d(y_i,l_i)\le 2n_0+r$ and
we obtain the contradiction with the asymptotic disjointness of
the sets $O(Z)$ and $V_\tau^{X_0}(L)$.

2. we can suppose that $z_i$, $s_i\notin \cup_{j=0}^{n_0}X_j$
where $k_n>r$ for each $n\ge n_0$. Then we have $z_i=s_i$ for each
$i\in\N$ and we can choose $y_i\in O(Z)\cap V_\tau^{X_0}(L)$.
Moreover, since $Z$ and $X\setminus
V_{\sigma_0}(X_0)\supset\cup_{j=1}^\infty K_j$ are asymptotically
disjoint, we can assume that $d(y_i,x_0)\to\infty$. We obtain the
contradiction again and the sets $Z$ and $S(V_\tau^{X_0}(L))$ are
asymptotically disjoint. Hence $x\notin(S(V_\tau^{X_0}(L)))'$ and
we have $\cap_{\tau\in\omega^\omega}(S(V_\tau^{X_0}(L)))'\cap
V_{\sigma_0}(X_0)'\subset (S_2(L))'$.

Put $V_A=\cup_{\tau\in\omega^\omega}\ex C_\tau$ and
$V_B=\cup_{\tau\in\omega^\omega}\ex D_\tau$.  Choose any point
$x\in V_A\cap V_B$. Then there exist $\tau_1$,
$\tau_2\in\omega^\omega$ such that $x\in\ex C_{\tau_1}\cap\ex
D_{\tau_2}$. Put $\tau=\max\{\tau_1,\tau_2\}$. Since $\ex
C_\tau\cap \ex D_\tau=\emptyset$, we have that $x\in(\ex
C_{\tau_1}\setminus \ex C_\tau)\cup(\ex D_{\tau_2}\setminus \ex
D_\tau)$. Consider the case when $x\in(\ex C_{\tau_1}\setminus \ex
C_\tau)$. Choose $n\in\N$ such that $A_{\tau_1}\setminus
A_\tau\subset B_n(x_0)$. Then $x\in S(B_n(x_0))'$. The same we
have in the case when $x\in(\ex D_{\tau_2}\setminus \ex D_\tau)$.
We obtain that $V_A\cap V_B\subset \cup_{j=1}^\infty
K_j'\subset\nu X\setminus V_{\sigma_0}(X_0)'$.

We have that $V_A\supset O_A\supset A'\cap X_0'$ and $V_B\supset
O_B\supset B'\cap X_0'$. Then there exists $\tau\in\omega^\omega$
such that $V_A\supset V_\tau(X_0)'\cap A'$ and $V_B\supset
V_\tau(X_0)'\cap B'$. Moreover, we can suppose that $(S_2(L)\cap
V_\tau(X_0))'\cap (A'\cup B')=\emptyset$ and $\sigma_0\le^*\tau$.
Then $U_A^1\cap U_B^2=\emptyset$ where $U_A^1=V_A\cap
V_\tau(X_0)'$, $U_B^1=V_B\cap V_\tau(X_0)'$. Then $U_A^1$ and
$U_B^1$ are two disjoint open subset of $\nu X$ such that
$U_A^1\supset V_\tau(X_0)'\cap A'$ and $U_B^1\supset
V_\tau(X_0)'\cap B'$. Since $\Ind(\nu X\setminus V_\tau(X_0)')\le
0$, there exists two open in $\nu X\setminus V_\tau(X_0)'$
disjoint sets $U_A^2$ and $U_B^2$ such that $U_A^2\supset (\nu
X\setminus V_\tau(X_0)')\cap A'$, $U_B^2\supset (\nu X\setminus
V_\tau(X_0)')\cap B'$ and $\nu X\setminus V_\tau(X_0)'=U_A^2\cup
U_B^2$. Since  $V_\tau(X_0)'$ is clopen, $U_A^2$ and $U_B^2$ are
open in $\nu X$.

The sets $U_A=U_A^1\cup U_A^2$ and $U_B=U_B^1\cup U_B^2$ are open
disjoint subsets of $\nu X$ such that $A'\subset U_A$, $B'\subset
U_B$ and $\nu X\setminus (U_A\cup U_B)\subset (S_2(L)\cap
V_\tau(X_0))'$.

Then $S_2(L)\cap V_\tau(X_0)$ is an asymptotic separator in $X$
between $A$ and $B$. Since $L$ is a kernel of $S_2(L)\cap
V_\tau(X_0)$, we have $\asInd(S_2(L)\cap V_\tau(X_0))<n$ by the
inductive assumption and the theorem is proved.
\enddemo

{\bf 3.} In this section we prove the subspace and addition
theorems.

\proclaim {Theorem 2} Let $X$ be a proper metric space and
$Y\subset X$. Then $\asInd Y\le\asInd X$.
\endproclaim

\demo{Proof} We shall apply induction with respect to $\asInd X$.
The case when $\asInd X=-1$ is trivial. The case when $\asInd X=0$
follows from the equivalence $\asInd X=0$ and $\Ind \nu X=0$ for
each metric proper space $X$ [2].

Let us assume that we have proved the theorem for each $X$ with
$\asInd X<n\ge 1$. Consider a proper metric space $X$ with $\asInd
X=n$. Suppose that $A$ and $B$ are asymptotically disjoint subsets
of $Y$ and $C$ is an asymptotic separator in $X$ between $V_A$ and
$V_B$ with $\asInd C<n$ where $V_A$ and $V_B$ are asymptotic
neighborhoods of $A$ and $B$ respectively such that $V_A$ and
$V_B$ are asymptotically disjoint. We will build an asymptotic
separator $L$ in $Y$ between $A$ and $B$ such that $\asInd L<n$
using a construction from [4, Lemma 5.4].

Denote $D_0=C$. Put $Z=Y\setminus(V_A\cup V_B)$. Let
$D_k=(N_k(D_0)\cap Z)\setminus B_k(\cup_{i=0}^{k-1}D_i)$  for
$k\in\N$. There exists a subset $L_k$ of $D_k$ such that $L_k$ is
$k$-discrete and for each $x\in D_k$ there exists $y\in L_k$ such
that $d(x,y)\le k$. Put $L=\cup_{i=1}^\infty L_i$. It was shown in
[4] that $L$ is an asymptotic separator in $Y$ between $A$ and
$B$. It is easy to see that $C$ is a kernel in $C\cup L$. Hence
$\asInd (C\cup L)<n$ and we have $\asInd L<n$ by the inductive
assumption. The theorem is proved.
\enddemo

\proclaim {Lemma 4} Let $X$ be a proper metric space and
$f:cX\to[0,1]$ is a continuous function. Consider any $a$,
$b\in[0,1]$ such that $a<b$. Then $(f^{-1}([a,b])\cap
X)'=(f^{-1}([0,b])\cap X)'\cap(f^{-1}([a,1])\cap X)'$.
\endproclaim

\demo{Proof} Choose any $x\in(f^{-1}([0,b])\cap
X)'\cap(f^{-1}([a,1])\cap X)'$ and any neighborhood $V$ of $x$ in
$cX$. Then there exist two sequences $(a_i)$, $(b_i)$ in $V\cap X$
such that $f(a_i)\ge a$, $f(b_i)\le b$ , $d(a_i,b_i)\le r$ for
some $r>0$ and $a_i$, $b_i\in X\setminus B_i(x_0)$ for each
$i\in\N$. We can suppose that $f(a_i)\to c_1$ and $f(b_i)\to c_2$.
Since $f|X$ is slowly oscillating, then $c_1=c_2=c\in[a,b]$. Then
we have $a<c\le b$ or $a\le c<b$. Consider the case $a\le c< b$.
There exists $n_0\in\N$ such that $f(b_n)\in[a,b]$ for each $n\ge
n_0$. Hence $\cl V\cap(f^{-1}([a,b])\cap X)'\neq\emptyset$ and
$x\in(f^{-1}([a,b])\cap X)'$. The proof is analogous in the case
$a<c\le b$. The inclusion $(f^{-1}([a,b])\cap
X)'\subset(f^{-1}([0,b])\cap X)'\cap(f^{-1}([a,1])\cap X)'$ is
trivial and the lemma is proved.
\enddemo

\proclaim {Theorem 3} Let $X$ be a proper metric space and
$X=Y\cup Z$ where $Y$ and $Z$ are unbounded sets. Then $\asInd
X\le\asInd Y+\asInd Z$.
\endproclaim

\demo{Proof} We shall apply induction with respect to $\asInd
Y+\asInd Z$. If $\asInd Y+\asInd Z=0$, then $\asInd Y=\asInd Z=0$.
We have that $\Ind Y'=\Ind Z'=0$ [2]. Since $\nu X=Y'\cup Z'$, we
have $\Ind\nu X=0$ [3, Theorem 2.2.7]. Then we have $\asInd X=0$
[2].

Assume that the theorem is proved for each $Y$ and $Z$ with
$\asInd Y\le m\ge 0$, $\asInd Z\le l\ge 0$ and $m+l<n\ge 1$.
Consider the case when $m+l=n$. Let $A$ and $B$ be asymptotically
disjoint subsets of $X$. Then $A'\cap B'=\emptyset$. Choose a
continuous function $f:\nu X\to [0,1]$ such that
$f(A')\subset\{0\}$ and $f(B')\subset\{1\}$ and extend it to a
continuous function $g:cX\to[0,1]$. Then the sets
$g^{-1}[0,\frac{1}{3}]\cap X$ and $g^{-1}[\frac{2}{3},1]\cap X$
are asymptotically disjoint. We can choose an asymptotic separator
$L_1$ in $Y$ between $g^{-1}[0,\frac{1}{3}]\cap Y$ and
$g^{-1}[\frac{2}{3},1]\cap Y$ such that $\asInd L_1<m$. Put
$L_2=g^{-1}[\frac{1}{3},\frac{2}{3}]\cap Z$. Then $\asInd
L_2\le\asInd Z\le l$. We have that $\asInd L_1\cup L_2\le m-1+l<n$
by the inductive assumption.

Let us show that $L_1\cup L_2$ is an asymptotic separator between
$A$ and $B$. It is easy to see that
$K=g^{-1}[\frac{1}{3},\frac{2}{3}]\cap X$ is an asymptotic
separator between $A$ and $B$ in $X$. Put
$V_1=\ex(g^{-1}[0,\frac{1}{3})\cap X)$ and
$V_2=\ex(g^{-1}(\frac{2}{3},1]\cap X)$. We have that $V_1$ and
$V_2$ are disjoint open sets in $\nu X$ such that $A'\subset V_1$
and $B'\subset U_1$. It follows from Lemma 4 that $\nu X\setminus
K'\subset V_1\cup U_1$. Since $L_1$ is an asymptotic separator in
$Y$ between $g^{-1}[0,\frac{1}{3}]\cap Y$ and
$g^{-1}[\frac{2}{3},1]\cap Y$, there exist two open in $Y'$
disjoint sets $O_1$ and $O_2$ such that
$O_1\supset(g^{-1}[0,\frac{1}{3}]\cap Y)'$,
$O_2\supset(g^{-1}[\frac{2}{3},1]\cap Y)'$ and $Y'\setminus
L_1'=O_1\cup O_2$. Put $V_2=O_1\setminus Z'$ and $U_2=O_2\setminus
Z'$. Then the sets $V=V_1\cup V_2$ and $U=U_1\cup U_2$ are open
disjoint subsets of $\nu X$ such that $A'\subset V$, $B'\subset U$
and $\nu X\setminus (L_1\cup L_2)'\subset V\cup U$. The theorem is
proved.
\enddemo

{\bf 4.} We investigate the transfinite extension of $\asInd$ in
this section. Let us recall the definition of the transfinite
 large inductive dimension $\trInd$ [3]: $\trInd X=-1$ iff $X=\emptyset$;
 $\trInd X\le \alpha$ for an ordinal number if for every two
disjoint closed subsets $A_1, A_2\subset X$ there is a partition
$C$ with $\Ind C\le \alpha$.

Define the transfinite extension $\trasInd X$ analogously:
$\trasInd X=-1$ if and only if $X$ is bounded; $\trasInd X\le
\alpha$ where $\alpha$ is an ordinal number if for every two
asymptotically disjoint sets $A$, $B\subset X$ there is an
asymptotic separator $C$ with $\trasInd C\le \beta$ for some
$\beta<\alpha$. Naturally we say $\trasInd X=\alpha$ if $\trasInd
X\le \alpha$ and it is not true that $\trasInd X\le \beta$ for
some $\beta<\alpha$. We set $\trasInd X=\infty$ if for each
ordinal number $\alpha$ it is not true that $\trasInd X\le
\alpha$.

The proof of the following theorem is the same as the one of
Theorem 2.

\proclaim {Theorem 4} Let $X$ be a proper metric space and
$Y\subset X$. Then $\trasInd Y\le\trasInd X$.
\endproclaim

Let $X$ be a metric proper space and $\{A_i|i\in\N\}$ is a
countable family of subsets. We say that the family
$\{A_i|i\in\N\}$ is {\it asymptotically discrete} if for each
$i\in\N$ the sets $A_i$ and $\cup_{j\neq i} A_j$ are
asymptotically disjoint. We say that a proper metric space $X$ is
{\it asymptotically $S$-like} if $X$ can be represented as the
union of a sequence $X_1, X_2,\dots$ of subsets such that $\asInd
X_i\ge i$ and the family $\{X_i|i\in\N\}$ is asymptotically
discrete. The class of all asymptotically $S$-like spaces we
denote by $\S$.

\proclaim {Lemma 5} If a proper metric space $X$ is asymptotically
$S$-like, then $\trasInd X=\infty$.
\endproclaim

\demo{Proof} Suppose the contrary. Then there exists $X\in\S$ such
that $\trasInd X<\infty$. Put $\xi=\min\{\alpha\mid$ there exists
$X_\alpha\in\S$ such that $\trasInd X_\alpha=\alpha\}$. Clearly,
$\xi\ge\omega$ where $\omega$ is the first infinite ordinal
number.

Choose any $X\in\S$ such that $\trasInd X=\xi$. Let us represent
$X$ as the union $\cup_{i=1}^\infty X_i$ of subsets $X_i$ such
that $\asInd X_i\ge i$ and the family $\{X_i|i\in\N\}$ is
asymptotically discrete. For each $i\in\N$ choose two
asymptotically disjoint subsets $A_i$ and $B_i$ of $X_i$ such that
for each asymptotic separator $L_i$ in $X_i$ between $A_i$ and
$B_i$ we have $\trasInd L_i\ge i-1$.

We build by induction sets $C_i\subset A_i$ and $D_i\subset B_i$
such that $d(C_i,\cup_{j=1}^i D_j)\ge i$ and $d(D_i,\cup_{j=1}^i
C_j)\ge i$ for each $i\in\N$. Since $A_1$ and $B_1$ are
asymptotically disjoint, there exists $r>0$ such that
$d(C_1,D_1)\ge 1$ where $C_1=A_1\setminus B_r(x_0)$ and
$D_1=B_1\setminus B_r(x_0)$. Suppose we have built $C_k$ and $D_k$
for each $k\le n\ge 1$. Build $C_{n+1}$ and $D_{n+1}$. Since the
sets $\cup_{j=1}^n X_j$ and $X_{n+1}$ are asymptotically disjoint,
there exists $r>0$ such that $d((\cup_{j=1}^n X_j)\setminus
B_r(x_0),X_{n+1}\setminus B_r(x_0))\ge n+1$. Since the sets
$A_{n+1}$ and $B_{n+1}$ are asymptotically disjoint, there exists
$t>0$ such that $d(A_{n+1}\setminus B_r(x_0),B_{n+1}\setminus
B_r(x_0))\ge n+1$. Put $s=\max\{r,t\}$ and
$C_{n+1}=A_{n+1}\setminus B_s(x_0)$, $D_{n+1}=B_{n+1}\setminus
B_s(x_0)$. It is easy to see that the sets $C=\cup_{i=1}^\infty
C_i$, $D=\cup_{i=1}^\infty D_i$ are asymptotically disjoint and
for each asymptotic separator $L_i$ in $X_i$ between $C_i$ and
$D_i$ we have $\trasInd L_i\ge i-1$.

Choose any asymptotic separator $L$ in $X$ between $C$ and $D$
such that $\trasInd L<\xi$. Since the family $\{X_i\}$ is
asymptotically discrete, we have that $L_i=L\cap X_i$ is an
asymptotic separator in $X_i$ between $C_i$ and $D_i$ for each
$i\in\N$. Then $\asInd L_i\ge i-1$. Thus, $L\in\S$ and we have a
contradiction. The lemma is proved.
\enddemo

\proclaim {Lemma 6} Let $X$ be a proper metric space such that
$\trasInd X<\infty$. Then for each $x\in\nu X$ there exists a
neighborhood $V$ of $x$ in $cX$ such that $\asInd V\cap X<\infty$.
\endproclaim

\demo{Proof} Suppose the contrary. Then there exists $x\in\nu X$
such that $\asInd V\cap X=\infty$ for each neighborhood $V$ of $x$
in $cX$. Let us build by induction the sequence $(L_i)$ of subsets
of $X$ and the sequence $(V_i)$ of neighborhoods of $x$ in $cX$
such that the sets $ L_k$ and $V_k\cap X$ are asymptotically
disjoint for each $k\in\N$, $L_{n+1}\subset V_n$, $\asInd
L_{n}=\infty$  and $V_1\supset\dots\supset V_n\supset
V_{n+1}\supset\dots$.

We have, particularly, that $\asInd X=\infty$. There exist two
asymptotically disjoint subsets $A$, $B\subset X$ such that
$\asInd L=\infty$ for each asymptotic separator $L$ between $A$
and $B$. We can assume that $x\notin B'$. Choose a continuous
function $f:\nu X\to [0,1]$ such that $f(x\cup A')\subset \{0\}$
and $f(B')\subset \{1\}$. Choose an asymptotic separator $L_1$
between the sets $f^{-1}[0,\frac{1}{3}]\cap X$ and
$f^{-1}[\frac{2}{3},1]\cap X$ and put
$V_1=f^{-1}[0,\frac{1}{3}]\cap X$.

Assume we have built $L_i$ and $V_i$ for each $i\le n\ge n+1$. We
have that $\trasInd V_n=\infty$. So, there exist two
asymptotically disjoint subsets $C$ and $D$ of $V_n$ such that
$\asInd L=\infty$ for each asymptotic separator $L$ between $C$
and $D$. We can choose $L_{n+1}$ and $V_{n+1}$ as before. The
sequences $(L_i)$  and  $(V_i)$ are built. It is easy to check
that the family ${L_i}$ is asymptotically discrete. Put
$L=\cup_{i=1}^\infty L_i$.  So, $\trasInd L=\infty$ by Lemma 5 and
$\trasInd X=\infty$ by Theorem 4.
\enddemo

\proclaim {Theorem 5} Let $X$ be a proper metric space such that
$\trasInd X<\infty$. Then $\asInd X<\infty$.
\endproclaim

\demo{Proof}We can choose for each point $x\in\nu X$ a
neighborhood $V_x$ in $cX$ such that $\asInd V_x\cap X<\infty$.
Since $\nu X$ is compact, we have $\nu X\subset\cup_{i=1}^k V_i=V$
where each $V_i$ is an open in $cX$ set with $\asInd V_i\cap
X<\infty$. Then $\asInd V\cap X<\infty$ by Theorem 3. Moreover,
there exists $r>0$ such that $X\subset V\cup B_r(x_0)$. Hence
$\asInd X<\infty$.
\enddemo

\enddocument